\documentclass{article}
\usepackage{amsfonts,xr,graphicx,amsmath}

\begin{document}
\title{{Sets of unique continuation for heat equation} }
 \author{{Nikolai Nadirashvili\thanks{ Aix-Marseille Universit\'{e}, I2M, 39, rue F. Joliot-Curie, 13453
Marseille, FRANCE, nicolas@cmi.univ-mrs.fr}, \hskip .4 cm Nadezda Varkentina\thanks{Universit\'{e} Bordeaux 1, 351, Cours de la Lib\'{e}ration, 33405 Talence Cedex FRANCE, nadezda.varkentina@u-bordeaux1.fr}}}

\date{}
\maketitle

\def\n{\hfill\break} \def\al{\alpha} \def\be{\beta} \def\ga{\gamma} \def\Ga{\Gamma}
\def\om{\omega} \def\Om{\Omega} \def\ka{\kappa} \def\lm{\lambda} \def\Lm{\Lambda}
\def\dl{\delta} \def\Dl{\Delta} \def\vph{\varphi} \def\vep{\varepsilon} \def\th{\theta}
\def\Th{\Theta} \def\vth{\vartheta} \def\sg{\sigma} \def\Sg{\Sigma}
\def\bendproof{$\hfill \blacksquare$} \def\wendproof{$\hfill \square$}
\def\holim{\mathop{\rm holim}} \def\span{{\rm span}} \def\mod{{\rm mod}}
\def\rank{{\rm rank}} \def\bsl{{\backslash}}
\def\il{\int\limits} \def\pt{{\partial}} \def\lra{{\longrightarrow}}

\def\C{\mathbb{C}}
\def\S{\mathbb{S}}
\def\Z{\mathbb{Z}}
\def\R{\mathbb{R}}
\def\N{\mathbb{N}}
\def\H{\mathbb{H}}
\def\tilde{\widetilde}
\def\epsilon{\varepsilon}

\def\n{\hfill\break} \def\al{\alpha} \def\be{\beta} \def\ga{\gamma} \def\Ga{\Gamma}
\def\om{\omega} \def\Om{\Omega} \def\ka{\kappa} \def\lm{\lambda} \def\Lm{\Lambda}
\def\dd{\delta} \def\Dl{\Delta} \def\vph{\varphi} \def\vep{\varepsilon} \def\th{\theta}
\def\Th{\Theta} \def\vth{\vartheta} \def\sg{\sigma} \def\Sg{\Sigma}
\def\bendproof{$\hfill \blacksquare$} \def\wendproof{$\hfill \square$}
\def\holim{\mathop{\rm holim}} \def\span{{\rm span}} \def\mod{{\rm mod}}
\def\rank{{\rm rank}} \def\bsl{{\backslash}}
\def\il{\int\limits} \def\pt{{\partial}} \def\lra{{\longrightarrow}}
\def\pa{\partial } 
\def\ra{\rightarrow }
\def\sm{\setminus }
\def\ss{\subset }
\def\ee{\epsilon }

 {\em Abstract.} We study  nodal lines of solutions to the heat equations. We are interested in
 the global geometry of nodal sets, in the whole domain of definition of the solution. The local
 structure of nodal sets is a well understander subject, while the global geometry of nodal
 lines is much less clear. We give a detailed analysis of a simple component of a nodal set of a
 solution of the heat equation. Our results are motivated by some applied problems. They
 related to classical problems of the unique continuation and the backward uniqueness for 
 parabolic equations.

  
\section{Introduction}
\bigskip

Let $u(x,t)$ be a bounded solution of the heat equation defined on a half-plane:  $\H= \R\times \R^+$,

\begin{equation}
\left\{
\begin{array}{ll}
u_{t} = u_{xx}\quad in\ \H,\\
u(x,0) = f(x),
\end{array}
\right.
\label{eq:1}
\end{equation}
where $f(x)$ is a continuous, bounded function. We are concerned with the following inverse problem.
Let  $s\subset \H$ be a subset of the half-plane. Assume that we know the restriction of $u $ on $s$. Can we recover then the function $f(x)$.  The question is motivated by different  applied problems, for instance it appears in interpretation of experimental data to defining the absorption of laser irradiation by a solid target~\cite{RefVarkentina}.
We will assume that $s$ is a curve in $\H$. Assume first that $s$ is a horizontal line, $\{t=1\}$.
If we assume that $u(x,t)$ is uniformly bounded, then $ 2 \sqrt \pi u(x,1)$ is given by the convolution:
$f(x)\ast \exp{-x^2/4}$
and the function $f(x)$ can be recovered by the operation of standard deconvolution of $u(x,1)$ with Gauss function.\\
Without assumption on boundness of u there is no uniqueness in determining of $f(x)$ as it follows from the classical Tychonov examples, see~\cite{RefJones, RefLittman}. However, it is sufficient for the uniqueness to assume the boundness of $u(x,t)$, for $x>0$,~\cite{RefEscauriaza}. Notice, that since $u(x,t)$ is an analytic function of $x$ it follows that for the uniqueness it is sufficient to take just a segment of the horizontal line for the class of the bounded solutions $u(x,t)$.

Let $ I = \lbrack 0,1\rbrack$. Then the following result holds, see, e.g., ~\cite{RefEvans}.

\emph{Proposition}. {\it Assume that for x $\in$ I, u(x,1)=0, then  $ f(x) \equiv$ 0.}
 
Results of such type are known as backward uniqueness for the heat equation.

In this note we are going to generalize \emph{Proposition} taking instead of the segment I more general set $\sigma$ on $(x,t)$-plane. We say that the set $\sigma$ is a set of  unique continuation for the heat equation if any bounded solution of equation (\ref{eq:1})  which vanishes on $\sigma$ vanishes also for $t=0$. Notice that any open set is always a unique continuation set for very general
parabolic equations, \cite{RefSogge} \\
We show that the graph of any bounded function of $x$ is a set of unique continuation.

Let $g(x), g>0$ be a continuous bounded function defined on $x\in \R^+$. Denote by $\widetilde{g} $ the graph of the function $g$. Is $\widetilde{g}$ the set of unique continuation?

We prove two following theorems.\\

\textbf{Theorem 1.} \label{th:1}
\emph{Assume that $g>0$ is a bounded continuous function on $(0, \infty)$. Let u be a solution of the heat equation (\ref{eq:1}) defined in $ \H$ and u vanishes on $\widetilde{g}$:
$$u(x, g(x)) = 0.$$
Then $ u(x,0) \equiv 0$.\\ }

 \emph{Theorem (\ref{th:1})} gives a description of the structure of nodal set of solutions to the heat equation. The local structure of the nodal set of solutions of parabolic equations was intensively studied, see,~\cite{RefLin, RefHanLin, RefAngenent, RefMatano}. The global geometry of nodal sets 
 is less understood then its local structure. For elliptic equations global geometry of nodal lines studied in ~\cite{RefEJN}. We prove a result on a global structure of the nodal line of the heat equation.  Let $u$ to be a solution of equation (\ref{eq:1}). Let $\Gamma \subset \H$ be the nodal set of $u$, i.e., $\Gamma$ is the set of zeros of $u$. Let $\gamma \subset \Gamma$ be a component of the nodal set. 
We say that $\gamma$ is a nodal curve if $\gamma$ is a simple curve. The following result together with \emph{Theorem (\ref{th:1})} gives a characterization of nodal curves to solutions of the heat equation.
\\

\textbf{Theorem 2.}
\emph{Let $\gamma \subset \H$ be a nodal curve of a bounded solution of the heat equation (\ref{eq:1}). Then $\gamma$ has an end point on the line $t=0$.}
\\

\section{Proof of the theorems}

We will use the Poisson representation of solutions of the problem (\ref{eq:1}), see, e.g.,~\cite{RefEvans}. 
Let $u$ be a solution of problem (\ref{eq:1}), then:

\begin{equation}
u(x, t) = \frac{1}{\sqrt {4\pi t}} \int_{-\infty}^{\infty}e^{-\frac{(x-y)^2}{4t}} f(y)dy
\label{eq:3}
\end{equation} 

We will use  the complexification of solutions of the heat equation,
i.e., we consider equation (\ref{eq:1}) over the complex values of the variable x. Formula (\ref{eq:3}),  give a ready extension of solutions of the equation over the complex field $\mathbb{C}$.\\

Substituting complex values x=z $\in\mathbb{C}$ into (\ref{eq:3}) and integrating over $\mathbb{R}$ we immediately get the following lemma.\\

\textbf{Lemma 1.}
\emph{For every $t>0$ function u(z,t) is defined for all z $\in\mathbb{C}$.
Moreover, $u(z,t)$ satisfies the following inequality:
$$|u(z,t)|\leq\frac{M}{2\pi t} e^{(\Im z)^2/4t}.$$}

Let $Q\subset \mathbb{R}^2$ be a square $\{\vert x_{1}\vert <1,  0< x_{2}<2\}$. Denote by 
$\gamma_i$, $i=1,\ldots,4$ the edges of $ Q$, such that $\gamma_1 $ is on the line $\{x_2=0\}$.\\

\textbf{Lemma 2.} \emph{ Let $v$ be a subharmonic function in $Q$. Assume that $v\leq a_i$
on $\gamma_i$, then:
$$v(0,1)\leq \frac 14 \sum a_i.$$} 

Proof.
Let $\widetilde {v}$ be an averaging of ${v}$ over the rotation of Q, then ${v}\leq \frac 14 \sum a_i$ on $\pa Q$. By the maximum principle $v(0)=\widetilde {v}\leq \frac 14 \sum a_i$.

Denote $l=\{x_1,x_2: x_1=0, 0<x_2<1\}$.

\textbf{Lemma 3.} \emph{ Let $v$ be a subharmonic function in $Q$. Assume that $v\leq s$
on $\gamma_1$ and $v\leq p$ on $\gamma_i,\, i=2,3,4$, where $s<p$, then
$$v\leq \frac14 (s+3p)\quad on \quad l.$$}

Proof. Denote by $Q^r\subset \mathbb{R}^2$ the square $\{\vert x_{1}\vert <r,  0< x_{2}<2r\}$, $0<r<1$. Since $v$ is a subharmonic function in $Q$ then by the maximum principle $v<p$ on $\pa Q^r$. From the inequality $v<s$ on $\pa Q^r\cap \gamma_1$ and from \textbf{Lemma 2} it follows that $v<(s+3p)/4$ at the centre of $Q^r$. Since the union of all centers of the squares
$Q^r$, $0<r<1$ coincide with the segment $l$ the lemma follows.

 We need the following 
Phragm\'en-Lindel\"of type lemma.

\textbf{Lemma 4.}
\emph{Let w be an entire holomorphic function on the complex plane $\mathbb{C} = \{z=x+iy\}$. Assume that $\vert w(z)\vert \leq Ce^{A\vert z\vert^2}$, where C, A are positive constants. Assume that for any $N>0$ there exists a constant $C_N>0$ such that for $ x>0$ $\vert w(x,0)\vert < C_Ne^{-Nx^2}$, then $w \equiv 0$.}\\

Proof. Assume by contradiction that $w$ is not identically zero on $\C$. Set $v=\ln |w|$, then $v$ is a subharmonic function on the plane (with poles at zeros of $w$).
Moreover, $v(z) <C+A|z|^2$ on $\C$ and for any $N>0$ there exists a constant $C_N>0$ such that
$v(x,0)<C_N-N|z|^2$, $x>0$.

Denote
$$Q_L=\{x,y: L<x<3L,\, 0<y<2L\},$$
$$q_L=\{x,y: L<x<3L, \, y=0\},$$
$$l_L=\{x,y: x=2L,\, 0<y<L.\}$$

Then from the last inequalities it follows that
$$v<C+13L^2\quad on \quad \pa Q_L,$$
$$v<-N(L)L^2\quad on \quad q_L,$$
where $N(L)\ra +\infty $ as $L\ra +\infty$. Thus, by \textbf{Lemma 3}:
$$v<-M(L)L^2\quad on \quad l_L,$$
where $M(L)\ra +\infty $ as $L\ra +\infty$. Similarly, the last inequality holds on the segment:
$$l_L^-=\{x,y: x=2L,\, -L<y<0.\}$$

Now we can apply \textbf{Lemma 3} to the squares: 
$$S_L=\{x,y: |x|<2L,\, |y|<2L\}.$$
We have:
$$v<-M(L)L^2\quad on \quad l_L\cap l_L^-,$$
$$v<C+4AL^2 \quad on \quad \pa S_L.$$
By \textbf{Lemma 3}
$$v(0)< \frac 14 (3C+12AL^2-M(L)L^2).$$
Taking sufficiently large $L>0$ we get that $v(0)$ is less than any given negative constant.
Since we may assume without loss that $v(0)$ is finite, we got a contradiction and the lemma follows.

As a corollary of \textbf{Lemma 1} and \textbf{Lemma 4},  we have \textbf{Lemma 5.}

\textbf{Lemma 5.}
\emph{Let u be a bounded solution of the heat equation (\ref{eq:1}). Assume that for any $N>0$ there is a constant $C_N>0$ such that for $ x>0 $ the following inequality holds: \\
$$\vert u(x,1)\vert < C_N e^{-Nx^2},$$
then u$\equiv$0.}

Proof of \textbf{Theorem 1.} First we prove that function $g$ can not have a local maxima or minima. Assume $g( c)=g(d),\, 0<c<d$. We prove that  for any $x_1\in (c,d),\, g( c ) =g(d)<g(x_1)$.  Assume by contradiction $g(d)\geq g(x_1)$. We may assume without loss that $g$ attains its
infimum on the interval $(c,d)$. If $g(x_1) = g(d)$ then the segment $[(c,g( c)),(d,g(d))]$ is
in the zero set of $u $ and thus since $u$ is a real analytic function in the $x$ variable $u$ vanishes
on the whole line $t=g( c )$. By the backward uniqueness, see [4],  $u\equiv 0$.   Therefore
$g(d)>g(x_1)$.   Denote by
$G\ss \H$ a domain, 
$$G= \{ (x,t): c<x<d, g(x)<t<g( c )\}.$$
Let $l$ be the segment $t=g( c ),\,  c\leq x\leq d$, $\gamma =\pa G \sm l$.
then $u=0$ on $\gamma$ and by the maximum principle for the heat equation, see, e.g., ~\cite{RefEvans}, $u=0$ in $G$. From unique continuation and from backward uniqueness
 for the heat equation it follows that $u\equiv 0$ at $\H$. 
 
Assume now that $g$ has a local maximum at a point $c$. Since $g$ has no local minimums
 $g$ is a monotone decreasing function on $(c, \infty )$. Denote:
 $$D=\{ (x,t): c<x,\, g(x)<t<g( c)\}.$$
 Since $u$ is vanishes on the graph of $g$ and $u$ is bounded in $\H$, then by the maximum principle $u\equiv 0 $ in $D$ (notice that the maximum principle holds also in unbounded domains,
 see, e.g., \cite{RefFJohn}. Hence $u\equiv 0$ in $\H$.

 By the assumptions of the theorem  $g$ is a bounded function on $\R$. 
Since $g$ has no local minima and maxima $g$ is a monotone function.  We may assume without loss that $g$ is a monotonically increasing function and $g$ tends to $1$ as $x$ goes to $+\infty$.  
Denote:
$$G_r=\{ (x,t): r<x, \, g(x)<t<1\}.$$
Denote by $P$ the parabolic boundary of $G_r$, i.e., $P$ is an intersection of $\pa G_r$ with the half-space $\{t<1\}$.

Without loss we may assume that $|u|<1$ in $\H$, then:
$$|u|<2\pi ee^{-\frac{(x-r)^2}{4(t-g( r ))}}\quad on \quad P,$$
and hence by the maximum principle the last inequality holds in $G_r$. In particular, on the 
line $\{t=1\}$ we have the following inequality:
$$|u(2r,1)|<2\pi ee^{-\frac{r^2}{4(1-g( r ))}}.$$
Since $g( r ) \to 1$ as $r\to +\infty$, the solution $u$ of the heat equation satisfies the assumptions of \textbf{Lemma 5}, thus it follows that $u\equiv 0$. The contradiction proofs the theorem.

Proof of \textbf{Theorem 2}. Assume by contradictn that $\gamma $ has no end points on the line 
$t=0$. Let $\Gamma(s): \R\to \gamma,\; \Gamma =(\Gamma_1,\Gamma_2)$ be a parametrisation of the nodal line $\gamma $. First we prove that $\gamma $ can be represented as a graph of a  function $g(x)$. Assume that there are $s_1,s_2\in \R, \, s_1<s_2$, such that $\Gamma_1(s_1)=\Gamma_2(s_2)$ and $0<\Gamma_1(s_2)<\Gamma_2(s_1)$. Then either
$\Gamma_1(s)\leq\Gamma_2(s_2)$ for all $s\geq s_2$, or there is $s_3>s_2$ such that
$\Gamma_1(s_3)>\Gamma_2(s_2)$. In the first case $u$ vanishes on the parabolic boundary of the domain $D$ bounded by
the line $t= \Gamma_2(s_3)$ and the corresponding part of $\gamma $ and hence $u\equiv 0$.
In the second case the same true for the domain bounded by the line $t= \Gamma_2(s_1)$
and a piece of $\Gamma (s),\, s\in [s_4,s_3],$ where $ s_4<s_3,\, \Gamma_2(s_4)=\Gamma_2(s_3)$.

Now as in the proof of \textbf{Theorem 1} we get that function $g$ has no local infimums. Thus $g$ can not be defined over a finite interval $(a,b)$, since in that case it tends to infinity at $a$ and $b$
and therefore should have a local minima on the interval. If $g$ is defined over semi-infinite interval
$(-\infty, a)$  or $(a, \infty)$  then since $g$ tends to infinity at $a$ it follows from the maximum pronciple that $u\equiv 0$.  For the case of  function $g$ defined over $\R$ the result follows from \textbf{Theorem 1}. 


\begin{thebibliography}{}
\bibitem{RefVarkentina}
N. Varkentina, {\it PhD thesis,} Aix-Marseille Universit\'e, 2012.
\bibitem{RefJones}
B. F. Jones, {\it A fundamental solution of the heat equation which is supported in a strip,} J. Math. Anal. Appl. \textbf{60}, (1977) 314-324
\bibitem{RefLittman}
V. Littman, {\it Boundary control theory for hyperbolic and parabolic partial differential equations with constant coefficients,} Annli Scuola Norm. Sup. Pisa \textbf{Serie IV, 3}, (1978) 567-580
\bibitem{RefEscauriaza}
L. Escauriaza, G. Seregin, V. $\mathrm{\breve{S}}$ver$\mathrm{\acute{a}}$k, {\it Backward uniqueness for parabolic equations,  }Arch. Ration. Mech. Anal. \textbf{169}, (2003) 147-157
%\bibitem{RefEvans}
%L. C. Evans, Joudkazis, V. Mizeikis, H. Mizawa,  J. of Applied Physics \textbf{106}, (2009) 051101-14
%\bibitem{RefDrude}
%N.W. Ashcroft, \textit{Solide State Physics} (Harcourt College Publishers, USA 1976) 2-27
\bibitem{RefEvans}
L. C. Evans, \textit{Partial Differential Equations} (AMS, Providence, 1998).

\bibitem{RefSogge}
C. D. Sogge, {\it A unique continuation theorem for second order parabolic differential operators,}
Ark. Mat. 28 (1990), 159Ð182.

\bibitem{RefLin}
F-H. Lin,  \textit{Nodal sets of solutions of elliptic and parabolic equations, } Comm. Pure Appl. Math., 44 (1991) 287Ð308

\bibitem{RefHanLin}
Q. Han, F-H. Lin,  
{\it Nodal sets of solutions of parabolic equations: II  } Comm. Pure Appl. Math., 
 47 (1994), 1219Ð1238

\bibitem{RefAngenent}
S.B. Angenent, 
{\it Nodal properties of solutions of parabolic equations}, Rocky Mountain J. Math. 21  (1991) 585-592

\bibitem{RefMatano}
H. Matano,
{\it Nonincrease of the lap number of solution for a one-dimensional semi-linear parabolic equation,}
J. Fac. Sci. Univ. Tokyo Sect. 1A Math. 29 (1982), 410-441

\bibitem{RefEJN}
A. Eremenko, D. Jakobson, N. Nadirashvili
{\it On nodal sets and nodal domains on $S^2$ and $R^2$}, Ann. Inst. Fourier, 57
(2007) , 2345-2360 


\bibitem{RefFJohn} F.John, {\it Partial Differential Equations,} 4th edition, Springer, 1982 


\end{thebibliography}
\end{document}